\documentclass[11pt,leqno]{amsart}
\usepackage{amssymb,verbatim,enumerate,ifthen}
\usepackage[mathscr]{eucal}
\usepackage[utf8]{inputenc}
\usepackage[T1]{fontenc}
\def\N{\mathbb{N}}
\def\R{\mathbb{R}}

\def\F{\mathscr{F}}

\def\M{\mathscr{M}}

\def\H{\mathscr{H}}
\def\E{{\mathscr{E}}}

\newtheorem{theorem}{Theorem}[section]
\newtheorem*{theorem*}{Theorem}
\def\Thm#1#2{\ifthenelse{\equal{#1}{*}}{\begin{theorem*}#2\end{theorem*}}
             {\begin{theorem}\label{T#1}#2\end{theorem}}}
\newtheorem{Atheorem}{Theorem}

\def\thm#1{Theorem~\ref{T#1}}
\newtheorem{proposition}[theorem]{Proposition}
\newtheorem*{proposition*}{Proposition}
\def\Prp#1#2{\ifthenelse{\equal{#1}{*}}{\begin{proposition*}#2\end{proposition*}}
{\begin{proposition}\label{P#1}#2\end{proposition}}}
\def\prp#1{Proposition~\ref{P#1}}

\newtheorem{corollary}[theorem]{Corollary}
\newtheorem*{corollary*}{Corollary}
\def\Cor#1#2{\ifthenelse{\equal{#1}{*}}{\begin{corollary*}#2\end{corollary*}}
             {\begin{corollary}\label{C#1}#2\end{corollary}}}
\def\cor#1{Corollary~\ref{C#1}}

\newtheorem{lemma}[theorem]{Lemma}
\newtheorem*{lemma*}{Lemma}
\def\Lem#1#2{\ifthenelse{\equal{#1}{*}}{\begin{lemma*}#2\end{lemma*}}
             {\begin{lemma}\label{L#1}#2\end{lemma}}}
\def\lem#1{Lemma~\ref{L#1}}

\theoremstyle{definition}
\newtheorem{remark}[theorem]{Remark}
\newtheorem*{remark*}{Remark}
\def\Rem#1#2{\ifthenelse{\equal{#1}{*}}{\begin{remark}\rm #2\end{remark}}
             {\begin{remark}\label{R#1}\rm #2\end{remark}}}

\newtheorem{example}[theorem]{Example}
\newtheorem*{example*}{Example}
\def\Exa#1#2{\ifthenelse{\equal{#1}{*}}{\begin{example*}\rm #2\end{example*}}
             {\begin{example}\label{Ex#1}\rm #2\end{example}}}

\def\eq#1{{\rm(\ref{E#1})}}
\def\Eq#1#2{\ifthenelse{\equal{#1}{*}}
  {\begin{equation*}\begin{aligned}#2\end{aligned}\end{equation*}}
  {\begin{equation}\begin{aligned}\label{E#1}#2\end{aligned}\end{equation}}}

\begin{document}
\begin{flushright}
\end{flushright}
\vspace{5mm}

\date{\today}

\title{On approximately monotone and approximately H\"older functions}

\author[A. R. Goswami]{Angshuman R. Goswami}
\address[A. R. Goswami]{Doctoral School of Mathematical and Computational Sciences, University of Debrecen, 
H-4002 Debrecen, Pf.\ 400, Hungary}
\author[Zs. P\'ales]{Zsolt P\'ales}
\address[Zs. P\'ales]{Institute of Mathematics, University of Debrecen, 
H-4002 Debrecen, Pf.\ 400, Hungary}
\email{\{pales,angshu\}@science.unideb.hu}


\thanks{The research of the first author was supported by the Hungarian Scientific Research Fund (OTKA) Grant K-111651 and by the EFOP-3.6.1-16-2016-00022 project. This project is co-financed by the European Union and the European Social Fund.}

\begin{abstract}
A real valued function $f$ defined on a real open interval $I$ is called $\Phi$-monotone if, for all $x,y\in I$ with $x\leq y$ it satisfies
$$
  f(x)\leq f(y)+\Phi(y-x),
$$
where $\Phi:[0,\ell(I)[\,\to\mathbb{R}_+$ is a given nonnegative error function, where $\ell(I)$ denotes the length of the interval $I$. If $f$ and $-f$ are simultaneously $\Phi$-monotone, then $f$ is said to be a $\Phi$-Hölder function.

In the main results of the paper, we describe structural properties of these function classes, determine the error function which is the most optimal one. We show that optimal error functions for $\Phi$-monotonicity and $\Phi$-Hölder property must be subadditive and absolutely subadditive, respectively. Then we offer a precise formula for the lower and upper $\Phi$-monotone and $\Phi$-Hölder envelopes. We also introduce a generalization of the classical notion of total variation and we prove an extension of the Jordan Decomposition Theorem known for functions of bounded total variations. 
\end{abstract}

\maketitle

\section{Introduction} 

The main concepts and results of this paper are distillated from the following elementary observations. Assume that $I$ is a nonempty interval and a function $f:I\to\R$ satisfies the following inequality 
\Eq{0}{
  f(x)\leq f(y)+\varepsilon (y-x)^p \qquad (x,y\in I,\,x<y)
}
for some nonnegative constant $\varepsilon$ and real constant $p\in\R$. That is, $f$ is nondecreasing with an error term described in terms of the $p$th power function. Clearly, if $\varepsilon=0$, then the above condition is equivalent to the nondecreasingness of $f$. Conversely, one can notice that every nondecreasing function $f$ satisfies \eq{0}. On the other hand, if $p=1$, then \eq{0} holds if and only if the function $g(x):=f(x)+\varepsilon x$ is nondecreasing and hence $f(x)=-\varepsilon x$ is a strictly decreasing solution of \eq{0}. If $p<1$, then the function $f(x):=-\varepsilon x^p$, $(x>0)$ is a strictly decreasing solution of inequality \eq{0} on the interval $I=\,]0,\infty[\,$. 

Surprisingly, for $p>1$, the situation is completely different. 
Fix $a<b$ in $I$, then choose $n\in\N$ arbitrarily, set $u:=(b-a)/n$ and apply inequality \eq{0} for the values $x:=a+(k-1)u$ and $y:=a+ku$. Then we get
\Eq{*}{
  f(a+(k-1)u)\leq f(a+ku)+\varepsilon u^p \qquad (k\in\{1,\dots,n\}).
}
Adding up these inequalities side by side for $k\in\{1,\dots,n\}$, after trivial simplifications, we arrive at
\Eq{*}{
  f(a)=f(a+0u)
  \leq f(a+nu)+n\varepsilon u^p 
  = f(b)+\varepsilon n^{1-p}(b-a)^p \qquad (n\in\N).
}
Upon taking the limit $n\to\infty$, it follows that
\Eq{*}{
  f(a)\leq f(b) \qquad(a,b\in I,\,a<b),
}
which shows that $f$ is nondecreasing. Therefore, for $p>1$ a function $f:I\to\R$ satisfies \eq{0} for some nonnegative $\varepsilon$ if and only if $f$ is nondecreasing.

Another motivation for our paper comes from the theory of approximate convexity
which has a rich literature, see for instance \cite{BorNag13}, \cite{BurHaz11}, \cite{BurHazJuh11}, \cite{DilHowRob99}, \cite{DilHowRob00},  \cite{DilHowRob02}, \cite{Ger88c}, \cite{GilGonNikPal17}, \cite{GonNikPalRoa14}, \cite{Gre52}, \cite{HyeUla52}, \cite{Haz05a}, \cite{Haz10}, \cite{HazPal04}, \cite{HazPal05}, \cite{JarLac09}, \cite{JarLac10}, \cite{Krz01}, \cite{Lac99}, \cite{Mak17}, \cite{MakHaz17}, \cite{MakPal10b}, \cite{MakPal11a}, \cite{MakPal12b}, \cite{MakPal12c}, \cite{MakPal12a}, \cite{MakPal13a}, \cite{MakPal13b}, \cite{Mro01}, \cite{MroTabTab08}, \cite{NgNik93}, \cite{NikPal01}, \cite{NikPal07}, \cite{Pal01c}, \cite{Rol79}, \cite{SpuTab12},  \cite{TabTabZol10b}, \cite{TabTabZol10a}.
In these papers several aspects of approximate convexity were investigated: stability problems, Bernstein--Doetsch-type theorems, Hermite--Hadamard type inequalities, etc. 

In the paper \cite{Pal03a}, the particular case $p=0$ of inequality \eq{0} was considered and the following result was proved: A function $f:I\to\R$ satisfies \eq{0} for some $\varepsilon\geq0$ with $p=0$ if and only if there exists a nondecreasing function $g:I\to\R$ such that $|f-g|\leq\varepsilon/2$ holds on $I$. In other words, certain approximately monotone functions can be approximated by nondecreasing functions.

The above described observations and results motivate the investigation of classes of functions that obey a more general approximate monotonicity and also the related Hölder property. In fact, the class of approximate Hölder functions was introduced in the paper \cite{MakPal12a}, but this property was only investigated in the related context of approximate convexity. In this paper, we describe structural properties of these function classes, determine the error function which is the most optimal one. We show that optimal error functions for approximate monotonicity and for the Hölder property must be subadditive and absolutely subadditive, respectively. Then we offer a precise formula for the lower and upper approximately monotone and Hölder envelopes and also obtain sandwich-type theorems. In last section, we introduce a generalization of the classical notion of total variation and we prove a generalization of the Jordan Decomposition Theorem known for functions of bounded variations.

\section{On $\Phi$-monotone and $\Phi$-Hölder functions}

Let $I$ be a nonempty open real interval throughout this paper and let $\ell(I)\in\,]0,\infty]$ denote its length. The symbols $\R$ and $\R_+$ denote the sets of real and nonnegative real numbers, respectively. 

The class of all functions $\Phi:[0,\ell(I)[\,\to\R_+$, called error functions, will be denoted by $\E(I)$. Obviously, $\E(I)$ is a convex cone, i.e., it is closed with respect to addition and multiplication by nonnegative scalars. In what follows, we are going to define four properties related to an error function $\Phi\in\E(I)$.

A function $f:I\to\R$ will be called \emph{$\Phi$-monotone} if, for all $x,y\in I$ with $x\leq y$, 
\Eq{H1}{
  f(x)\leq f(y)+\Phi(y-x).
}
If this inequality is satisfied with the identically zero error function $\Phi$, then we say that $f$ is \emph{monotone (increasing)}.
The class of all $\Phi$-monotone functions on $I$ will be denoted by $\M_\Phi(I)$.
We also consider the class of all functions that are $\Phi$-monotone for some error function $\Phi\in\E(I)$:
\Eq{*}{
  \M(I):=\bigcup_{\Phi\in\E(I)}\M_\Phi(I).
}

A function $f:I\to\R$ will be called \emph{$\Phi$-H\"older} if, for all $x,y\in I$, 
\Eq{H2}{
  |f(x)-f(y)|\leq\Phi(|x-y|).
}
The class of all $\Phi$-Hölder functions on $I$ will be denoted by $\H_\Phi(I)$.
The family of all functions that are $\Phi$-Hölder for some error function $\Phi\in\E(I)$ will be denoted by $\H(I)$:
\Eq{*}{
  \H(I):=\bigcup_{\Phi\in\E(I)}\H_\Phi(I).
}

\Prp{MC1}{Let $\Phi_1,\dots,\Phi_n\in\E(I)$ and $\alpha_1,\dots,\alpha_n\in\R_+$. Then
\Eq{*}{
  \alpha_1\M_{\Phi_1}(I)+\cdots+\alpha_n\M_{\Phi_n}(I)
  &\subseteq \M_{\alpha_1\Phi_1+\cdots+\alpha_n\Phi_n}(I).
}
In particular, for all functions $\Phi\in\E(I)$, the class $\M_\Phi(I)$ is convex. Furthermore, $\M(I)$ is a convex cone.}

\begin{proof}
To prove the first inclusion, let $f\in \alpha_1\M_{\Phi_1}(I)+\cdots+\alpha_n\M_{\Phi_n}(I)$. Then, there exist $f_1,\dots ,f_n$ belonging to $\M_{\Phi_1}(I),\dots,\M_{\Phi_n}(I)$, respectively, such that
\Eq{999}{
f=\alpha_1 f_1+\dots+\alpha_n f_n.
}
Then, for all $x,y\in I$ with $x\leq y$, we have
\Eq{*}{     
 f_i(x)\leq f_i(y)+\Phi_i(y-x) \qquad (i\in\{1,\dots,n\}).
}
Multiplying this inequality by $\alpha_i$ and summing up side by side, we will arrive at
\Eq{*}{     
 f(x)=\sum_{i=1}^{n} \alpha_i f_i(x) 
 \leq \sum_{i=1}^{n} \alpha_if_i(y)+\sum_{i=1}^{n} \alpha_i\Phi_i(y-x)
 = f(y)+\Phi(y-x),
}
where $\Phi:=\sum_{i=1}^{n}\alpha_i\Phi_i$.
This shows that $f\in\M_{\Phi}(I)$, which proves statement.

The additional statements are immediate consequences of what we have proved.
\end{proof}

The following result is the counterpart of the previous statement.

\Prp{HA1}{Let $\Phi_1,\dots,\Phi_n\in\E(I)$ and $\alpha_1,\dots,\alpha_n\in\R$. Then
\Eq{*}{
  \alpha_1\H_{\Phi_1}(I)+\cdots+\alpha_n\H_{\Phi_n}(I)
  &\subseteq \H_{|\alpha_1|\Phi_1+\cdots+|\alpha_n|\Phi_n}(I).
}
In particular, for all functions $\Phi\in\E(I)$, the class $\H_\Phi(I)$is convex and central symmetric, i.e., $\H_\Phi(I)$ is closed with respect to multiplication by $(-1)$. Furthermore, $\H(I)$ is a linear space.}

\begin{proof}
To prove the first inclusion, let $f\in \alpha_1\H_{\Phi_1}(I)+\cdots+\alpha_n\H_{\Phi_n}(I)$. Then, there exist $f_1,\dots ,f_n$ belonging to $\H_{\Phi_1}(I),\dots,\H_{\Phi_n}(I)$, respectively, such that \eq{999} holds.
Then, for all $x,y\in I$, we have
\Eq{*}{     
 |f_i(x)-f_i(y)|\leq\Phi_i(y-x) \qquad (i\in\{1,\dots,n\}).
}
Multiplying this inequality by $|\alpha_i|$ and summing up side by side, we will arrive at
\Eq{*}{     
 |f(x)-f(y)|=\bigg|\sum_{i=1}^{n} \alpha_i (f_i(x)-f_i(y))\bigg| 
 \leq \sum_{i=1}^{n} |\alpha_i|\!\cdot\!|f_i(x)-f_i(y)|
 \leq \sum_{i=1}^{n} |\alpha_i|\Phi_i(y-x)
 =\Phi(y-x),
}
where $\Phi:=\sum_{i=1}^{n}|\alpha_i|\Phi_i$.
This shows that $f\in\H_{\Phi}(I)$, which proves the statement.

The additional statements are immediate consequences of what we have proved.
\end{proof}

\Prp{H0}{Let $\Phi\in\E(I)$. Then
\Eq{H0}{
\H_\Phi(I)=\M_\Phi(I)\cap(-\M_\Phi(I)).
}
Furthermore,
\Eq{H0+}{
\H(I)=\M(I)\cap(-\M(I)).
}}

\begin{proof}
Assume that $f$ is a $\Phi$-H\"older function. Then, for any $x,y\in I$ with $x\leq y$, $f$  will satisfy the inequality \eq{H2} and hence the inequalities
\Eq{Mon}
{f(x)-f(y)\leq\Phi(y-x),\qquad f(y)-f(x)\leq\Phi(y-x).
}
Rearranging theses two inequalities, we have that both $f$ and $-f$ are $\Phi$-monotone. That is $f\in \M_\Phi(I)\cap(-\M_\Phi(I)).$

To show the inverse inclusion, let $f\in \M_\Phi(I)\cap(-\M_\Phi(I)).$  Due to the property of $\Phi$-monotonicity of the two classes of function, $f$ will satisfy the two inequalities in \eq{Mon}. Hence, inequality \eq{H2} holds for $x\leq y$. This inequality being symmetric in $x$ and $y$, we get that \eq{H2} is satisfied for all $x,y\in I$.

To verify \eq{H0+}, let first $f$ be a member of $\H(I)$. Then there exists $\Phi\in\E(I)$ such that $f\in\H_\Phi(I)$. In view of the first part, this implies that 
\Eq{*}{
  f\in\M_\Phi(I)\cap(-\M_\Phi(I))\subseteq\M(I)\cap(-\M(I)). 
}
Thus, we have shown the inclusion $\subseteq$ for \eq{H0+}.

For the reversed inclusion, let $f\in\M(I)\cap(-\M(I))$. Then there exist $\Phi_1,\Phi_2\in\E(I)$ such that $f\in\M_{\Phi_1}(I)$ and $-f\in\M_{\Phi_2}(I)$. Define $\Phi:=\max(\Phi_1,\Phi_2)$. Then, obviously, $f\in\M_{\Phi}(I)$ and $-f\in\M_{\Phi}(I)$, therefore, $f\in\H_\Phi(I)\subseteq\H(I)$. This completes the proof.
\end{proof}

We say that a family $\F$ of real valued functions is \emph{closed with respect to the pointwise supremum} if $\{f_\gamma:I\to\R\mid\gamma\in\Gamma\}$ is a subfamily of $\F$ with a pointwise supremum $f:I\to\R$, i.e., 
\Eq{sup}{
  f(x)=\sup_{\gamma\in\Gamma} f_\gamma(x)\qquad(x\in I),
}
then $f\in\F$. Similarly, we can define that a family $\F$ of real valued functions is \emph{closed with respect to the pointwise infimum}. A family $\{f_\gamma:I\to\R\mid\gamma\in\Gamma\}$ is called a \emph{chain} if, for all $\alpha,\beta\in\Gamma$, either $f_\alpha\leq f_\beta$ or $f_\beta\leq f_\alpha$ holds on $I$. We say that a family $\F$ of real valued functions is \emph{closed with respect to the pointwise chain supremum (chain infimum)} if $\{f_\gamma:I\to\R\mid\gamma\in\Gamma\}\subseteq\F$ is a chain with a pointwise supremum (infimum) $f:I\to\R$, then $f\in\F$.

\Prp{MC2}{Let $\Phi\in\E(I)$. Then the class $\M_\Phi(I)$ is closed under pointwise infimum and supremum. Furthermore, $\M_\Phi(I)$ is closed with respect to the pointwise liminf and limsup operations.}

\begin{proof}
Assume that $\{f_\gamma\mid\gamma\in\Gamma\}$ is a family of $\Phi$-monotone functions with a pointwise supremum $f:I\to\R$, i.e., \eq{sup} holds.
Let $x,y\in I$ be arbitrary with $x\leq y$. Then, by the $\Phi$-monotonicity property, for all $\gamma\in\Gamma$, we have that
\Eq{*}{
  f_\gamma(x)\leq f_\gamma(y)+\Phi(y-x)\leq f(y)+\Phi(y-x).
}
Taking the supremum of the left hand side with respect to $\gamma\in\Gamma$, we get
\Eq{*}{
  f(x)\leq f(y)+\Phi(y-x),
}
which shows that $f$ is $\Phi$-monotone. The proof of the assertion related to the pointwise infimum is similar, therefore it is omitted.

To obtain the statements with respect to the liminf and limsup operations, let $f:I\to\R$ be the upper limit of a sequence $f_n:I\to\R$. Then
\Eq{*}{
  f=\inf_{n\in\N} g_n, \qquad\mbox{where}\qquad g_n:=\sup_{k\geq n} f_k.
}
If all the functions $f_n$ are $\Phi$-monotone, then for all $n\in\N$, the function $g_n$ is $\Phi$-monotone. On the other hand, the sequence $(g_n)$ is decreasing, therefore $f$ is the pointwise chain infimum of $\{g_n\mid n\in\N\}$, thus $f$ is also $\Phi$-monotone. 

In a similar way, one can prove that the class of $\Phi$-monotone functions is closed under the liminf operation.
\end{proof}

As an immediate consequence, we can see that if $f$ is $\Phi$-monotone, then 
$f_+=\max(f,0)$ is also $\Phi$-monotone.

\Prp{H2}{Let $\Phi\in\E(I)$. Then the class $\H_\Phi(I)$ is closed under pointwise infimum and pointwise supremum. Consequently, $\H_\Phi(I)$ is closed with respect to the pointwise liminf and limsup operations.}

\begin{proof} 
Assume that $f:I\to\R$ is the pointwise supremum of a family $\{f_\gamma\mid\gamma\in\Gamma\}\subseteq\H_\Phi(I)$. By \prp{H0}, we have that $\pm f_\gamma\in\M_\Phi(I)$ holds for all $\gamma\in\Gamma$. In view of \prp{MC2}, this implies that
\Eq{*}{
   f=\sup_{\gamma\in\Gamma} f_\gamma\in\M_{\Phi}(I) \qquad\mbox{and}\qquad
   -f=\inf_{\gamma\in\Gamma} (-f_\gamma)\in\M_{\Phi}(I).
}
Therefore, $f\in \M_{\Phi}(I)\cap(-\M_{\Phi}(I))=\H_{\Phi}(I)$. The proof of the statement for the pointwise infimum is analogous.

The statements concerning liminf and limsup operations follow from the first part exactly in the same way as in the proof of \prp{MC2}.
\end{proof}

As a trivial corollary, we obtain that if $f$ is $\Phi$-Hölder, then 
$|f|=\max(f,-f)$ is also $\Phi$-Hölder.

\section{Optimality of the error functions}

In what follows, a function $\Phi\in\E(I)$ will be called \emph{subadditive} if, for all $u,v\in\R_+$ with $u+v<\ell(I)$, the inequality
\Eq{sub}{
  \Phi(u+v)\leq \Phi(u)+\Phi(v)
}
holds. Obviously, a decreasing function $\Phi\in\E(I)$ is automatically subadditive. Indeed, if $u,v\geq0$ with $u+v<\ell(I)$, then $u\leq u+v$ implies $\Phi(u+v)\leq\Phi(u)$, which, together with the nonnegativity of $\Phi(v)$, yields \eq{sub}. A stronger property of a function $\Phi\in\E(I)$ is the \emph{absolute subadditivity}, which is defined as follows: for all $u,v\in\R$ with $|u|,|v|,|u+v|<\ell(I)$, the inequality
\Eq{ASF}{
  \Phi(|u+v|)\leq \Phi(|u|)+\Phi(|v|)
}
is satisfied. It is clear that absolutely subadditive functions are automatically subadditive. On the other hand, we have the following statement.

\Lem{ASF}{If $\Phi\in\E(I)$ is increasing and subadditive, then it is absolutely subadditive.}

\begin{proof}
Le $u,v\in\R$ with $|u|,|v|,|u+v|<\ell(I)$. If $uv\geq0$, then $|u|+|v|=|u+v|<\ell(I)$ and the subadditivity implies
\Eq{*}{
   \Phi(|u+v|)=\Phi(|u|+|v|)\leq \Phi(|u|)+\Phi(|v|).
}
In the case $uv<0$, one can easily check that $|u+v|\leq\max(|u|,|v|)$. Therefore, the monotonicity property of $\Phi$ yields
\Eq{*}{
  \Phi(|u+v|)\leq\Phi(\max(|u|,|v|))
  =\max(\Phi(|u|),\Phi(|v|))\leq \Phi(|u|)+\Phi(|v|).
}
Thus, we have proved \eq{ASF} in both cases.
\end{proof}

We mention here another related notion, the concept of increasing subadditivity which was introduced in \cite{MakPal12a}: $\Phi\in\E(I)$ is called \emph{increasingly subadditive} if, for all $u,v,w\in[0,\ell(I)[\,$ with $u\leq v+w$ the inequality
\Eq{*}{
  \Phi(u)\leq\Phi(v)+\Phi(w)
}
holds. One can easily see that this property implies absolute subadditivity, but the converse is not true in general.

The simplest but important error functions are of the form 
\Eq{*}{
 \Phi_p(0):=0,\qquad \Phi_p(u):=u^p \qquad(u>0),
}
where $p\in\R$. Their subadditivity and absolute subadditivity is characterized by the following statement.

\Prp{p}{Let $p\in\R$. Then 
\begin{enumerate}[(i)]
 \item $\Phi_p$ is subadditive on $\R_+$ if and only if $p\in\,]-\infty,1]$.
 \item $\Phi_p$ is absolutely subadditive on $\R_+$ if and only if $p\in[0,1]$.
\end{enumerate}
}

\begin{proof} Let $p\leq1$. To show that $\Phi_p$ is subadditive, it is enough to check \eq{sub} for $\Phi=\Phi_p$ in the case $uv\neq0$. Then
\Eq{*}{
  1=\frac{u}{u+v}+\frac{v}{u+v}
  \leq \Big(\frac{u}{u+v}\Big)^p+\Big(\frac{v}{u+v}\Big)^p.
}
Multiplying this inequality side by side by $(u+v)^p$, we get
\Eq{p}{
  \Phi_p(u+v)=(u+v)^p\leq u^p+v^p=\Phi_p(u)+\Phi_p(v),
}
which shows the subadditivity of $\Phi_p$.

If $p>1$, then with $u:=v>0$ in \eq{p}, we get
\Eq{*}{
  \Phi_p(u+u)=(u+u)^p=(2u)^p=2^pu^p>2u^p=\Phi_p(u)+\Phi_p(u),
}
therefore, $\Phi_p$ cannot be subadditive (in fact, one can see that $\Phi_p$ is superadditive).

If $p\in[0,1]$, then $\Phi_p$ is increasing and also subadditive on $\R_+$ (by the first assertion), hence, by \lem{ASF}, it is also absolutely subadditive on $\R_+$.

If $p>1$, then $\Phi_p$ is not subadditive, hence, it is also not absolutely subadditive. If $p<0$, then with $u:=n+1$ and $v:=-n$, the absolute subadditivity of $\Phi_p$ would imply
\Eq{*}{
  1=\Phi_p(1)=\Phi_p(|u+v|)\leq \Phi_p(|u|)+\Phi_p(|v|)=(n+1)^p+n^p
}
for all $n\in\N$. Upon taking the limit $n\to\infty$ and using $p<0$, we arrive at the contradiction $1\leq 0$. Hence $\Phi_p$ cannot be absolutely subadditive.
\end{proof}

It is also not difficult to see that the class of subadditive functions as well as the class of absolutely subadditive functions is closed with respect to pointwise supremum. Therefore, for any $\Phi\in\E(I)$, there exists a largest subadditive function $\Phi^\sigma\in\E(I)$ and a largest absolutely subadditive function $\Phi^\alpha\in\E(I)$ which satisfy the inequality $\Phi^\sigma\leq\Phi$ and $\Phi^\alpha\leq\Phi$ on $[0,\ell(I)[\,$, respectively. The functions $\Phi^\sigma$ and $\Phi^\alpha$ will be called the \emph{subadditive envelope} (or \emph{subadditive minorant}) and the \emph{absolutely subadditive envelope} (or \emph{absolutely subadditive minorant}) of the function $\Phi$, respectively. Obviously, the equalities $\Phi=\Phi^\sigma$ and $\Phi=\Phi^\alpha$ are valid if and only if $\Phi$ is subadditive and absolutely subadditive, respectively. More generally, the functions $\Phi^\sigma$ and $\Phi^\alpha$ can be constructed explicitly from $\Phi$ by the following results. 

\Prp{2}{Let $\Phi\in\E(I)$ be an arbitrary function. Define the function $\Phi^\sigma:[0,\ell(I)[\,\to\R_+$ by
\Eq{*}{
\Phi^\sigma(u):=\inf\big\{\Phi(u_1)+\dots+\Phi(u_n)\mid
  n\in\N, u_1,\dots,u_n\in\R_+\colon u_1+\dots+u_n=u\big\}.
}
Then $\Phi^\sigma$ is the largest subadditive function which satisfies the inequality $\Phi^\sigma\leq\Phi$ on $[0,\ell(I)[\,$. Furthermore, $\Phi^\sigma(0)=\Phi(0)$ and if, additionally, $\Phi$ is increasing, then $\Phi^\sigma$ is also increasing, and hence $\Phi^\sigma=\Phi^\alpha$.}

\begin{proof}
First we are going to prove the subadditivity of $\Phi^\sigma$. Let $u,v\in\R_+$ such that $u+v\in[0,\ell(I)[\,$. Let $\varepsilon>0$ be arbitrary. Then there exist $n,m \in\N$ and $u_1,\dots,u_n,v_1,\dots,v_m\in\R_+$ such that 
\Eq{*}
{u=\sum_{i=1}^{n} u_i,\qquad v=\sum_{j=1}^{m} v_j, \qquad
\sum_{i=1}^{n}\Phi(u_i)<\Phi^\sigma(u)+\frac{\varepsilon}{2} 
\qquad\text{and}\qquad
\sum_{j=1}^{m}\Phi(v_j)<\Phi^\sigma(v)+\frac{\varepsilon}{2}.
}
We have that $u+v=\sum_{i=1}^{n} u_i+\sum_{j=1}^{m} v_j$. Therefore, by the definition of $\Phi^\sigma$ and by the last two inequalities, we get
\Eq{*}{
  \Phi^\sigma(u+v)\leq \sum_{i=1}^{n}\Phi(u_i) + \sum_{j=1}^{m}\Phi(v_j)
  <\Phi^\sigma(u)+\Phi^\sigma(v)+\varepsilon.
}
Since $\varepsilon$ is an arbitrary positive number, we conclude that $\Phi^\sigma(u+v)\leq \Phi^\sigma(u)+\Phi^\sigma(v)$, which completes the proof of the subadditivity of $\Phi^\sigma$. By taking $n=1$, $u_1=u$ in the definition of $\Phi^\sigma(u)$, we can see that $\Phi^\sigma(u)\leq\Phi(u)$ also holds for all $u\in[0,\ell(I)[\,$.

By the definition of $\Phi^\sigma$ and the inequality $\Phi(0)\geq0$, we have that \Eq{*}{
  \Phi^\sigma(0)=\inf_{n\in\N} n\Phi(0)=\Phi(0).
}

Now assume that $\Psi:[0,\ell(I)[\,\to\R_+$ is a subadditive function such that $\Psi\leq\Phi$ holds on $[0,\ell(I)[\,$. To show that $\Psi\leq\Phi^\sigma$, let $u\in[0,\ell(I)[\,$ and $\varepsilon>0$ be arbitrary. Then 
there exist $n \in\N$ and $u_1,\dots,u_n\in\R_+$ such that 
\Eq{v}
{u=u_1+\dots+u_n \qquad\text{and}\qquad
\Phi(u_1)+\dots+\Phi(u_n)<\Phi^\sigma(u)+\varepsilon. 
}
Then, due to the subadditivity of $\Psi$,
\Eq{*}{
  \Psi(u)\leq \Psi(u_1)+\dots+\Psi(u_n)
  \leq \Phi(u_1)+\dots+\Phi(u_n)<\Phi^\sigma(u)+\varepsilon.
}
By the arbitrariness of $\varepsilon>0$, the inequality $\Psi(u)\leq\Phi^\sigma(u)$ follows for all $u\in[0,\ell(I)[\,$, which was to be proved.

To verify the last assertion, let $u,v\in [0,\ell(I)[\,$ with $v<u$. Let $\varepsilon>0$ be arbitrary. Then there exist $n \in\N$ and $u_1,\dots,u_n\in\R_+$ such that \eq{v} is satisfied. Define $v_i:=\frac{v}{u}u_i$. Then $v_i\leq u_i$, hence $\Phi(v_i)\leq\Phi(u_i)$. On the other hand, $v_1+\dots+v_n=\frac{v}{u}(u_1+\dots+u_n)=v$, which implies that
\Eq{*}{
  \Phi^\sigma(v)
  \leq \Phi(v_1)+\dots+\Phi(v_n) 
  \leq \Phi(u_1)+\dots+\Phi(u_n)
  <\Phi^\sigma(u)+\varepsilon.
}
Passing the limit $\varepsilon\to0$, we arrive at the inequality $\Phi^\sigma(v)\leq\Phi^\sigma(u)$, which proves the increasingness of $\Phi^\sigma$. If this is the case, then $\Phi^\sigma$ is also absolutely subadditive, therefore, $\Phi^\sigma=\Phi^\alpha$.
\end{proof}

The following lemma is instrumental for the construction of the absolute convex envelope of a given error function.

\Lem{2-}{Assume that $\Phi\in\E(I)$ is absolutely subadditive. Then, for all $n\in\N$,
$u_1,\dots,u_n\in\R$ with $|u_1|,\dots,|u_n|,|u_1+\dots+u_n|<\ell(I)$,
\Eq{*}{
  \Phi(|u_1+\dots+u_n|)\leq \Phi(|u_1|)+\dots+\Phi(|u_n|).
}} 

\begin{proof} The statement is trivial for $n=1$. 
For $n\geq2$, we prove the assertion by induction. If $n=2$, then it is equivalent to the absolute subadditivity of $\Phi$. Let $n\geq2$ and assume that the statement holds for $n$ variables. Let $u_1,\dots, u_{n+1}\in\,]-\ell(I), \ell(I)[\,$ such that 
$u_0:=u_1+\dots+u_{n+1}\in\,]-\ell(I), \ell(I)[\,$. We may assume that $u_0\geq0$, (otherwise we can replace $u_i$ by $-u_i$ for all $i\in\{1,\dots,n+1\}$ in the argument). Then, for at least one $i\in\{1,\dots,n+1\}$, we have that $u_i\geq0$. By the symmetry, we may also assume that  $u_{n+1}\geq0$. Then, using the inequalities
$u_0\geq0$, $u_{n+1}<\ell(I)$, $u_{n+1}\geq0$ and $u_0<\ell(I)$, respectively, we get
\Eq{*}{
  -\ell(I)&\leq u_0-\ell(I)=u_1+\dots+u_n+(u_{n+1}-\ell(I))\\
  &<u_1+\dots+u_n\leq u_1+\dots+u_n+u_{n+1}=u_0<\ell(I).
}
This shows that $u_1+\dots+u_n\in\,]-\ell(I), \ell(I)[\,$. Now, applying the inequality \eq{ASF} with $u:=u_1+\dots+u_n$ and $v:=u_{n+1}$ and then the inductive assumption, it follows that
\Eq{*}{
  \Phi(|u_1+\dots+u_{n+1}|)\leq \Phi(|u_1+\dots+u_n|)+\Phi(|u_{n+1}|)
  \leq \Phi(|u_1|)+\dots+\Phi(|u_n|)+\Phi(|u_{n+1}|).
}
Therefore, the statement is valid also for $n+1$ variables.
\end{proof}

\Prp{2+}{Let $\Phi\in\E(I)$ be an arbitrary function. Define the function $\Phi^\alpha:[0,\ell(I)[\,\to\R_+$ by
\Eq{*}{
\Phi^\alpha(u):=\inf\big\{\Phi(|u_1|)+\dots+\Phi(|u_n|)\mid
  n\in\N, |u_1|,\dots,|u_n|<\ell(I),\, u_1+\dots+u_n=u\big\}.
}
Then $\Phi^\alpha$ is the largest absolutely subadditive function which satisfies the inequality $\Phi^\alpha\leq\Phi$ and hence $\Phi^\alpha\leq\Phi^\sigma$ on $[0,\ell(I)[\,$.}

\begin{proof}
First we are going to prove the absolute subadditivity of $\Phi^\alpha$. Let $u,v\in\R$ such that $|u|,|v|,|u+v|<\ell(I)$. Without loss of generality, we may assume that $u+v$ is nonnegative (otherwise, we replace $u$ and $v$ by $(-u)$ and $(-v)$ in the argument below). Let $\varepsilon>0$ be arbitrary. Then there exist $n,m \in\N$ and real numbers $u_1,\dots,u_n,v_1,\dots,v_m\in\,]-\ell(I),\ell(I)[\,$ such that 
\Eq{*}
{u=\sum_{i=1}^{n} u_i,\qquad v=\sum_{j=1}^{m} v_j, \qquad
\sum_{i=1}^{n}\Phi(|u_i|)<\Phi^\alpha(|u|)+\frac{\varepsilon}{2} 
\qquad\text{and}\qquad
\sum_{j=1}^{m}\Phi(|v_j|)<\Phi^\alpha(|v|)+\frac{\varepsilon}{2}.
}
We have that $u+v=\sum_{i=1}^{n} u_i+\sum_{j=1}^{m} v_j$. Therefore, by the definition of $\Phi^\alpha$ and by the last two inequalities, we get
\Eq{*}{
  \Phi^\alpha(u+v)\leq \sum_{i=1}^{n}\Phi(|u_i|) + \sum_{j=1}^{m}\Phi(|v_j|)
  <\Phi^\alpha(|u|)+\Phi^\alpha(|v|)+\varepsilon.
}
Since $\varepsilon$ is an arbitrary positive number, we conclude that $\Phi^\alpha(u+v)\leq \Phi^\alpha(|u|)+\Phi^\alpha(|v|)$, which completes the proof of the absolute subadditivity of $\Phi^\alpha$. By taking $n=1$ and $u_1=u$ in the right hand side expression for $\Phi^\alpha$, we can see that $\Phi^\alpha\leq\Phi$ holds.

Now assume that $\Psi\in\E(I)$ is an absolutely subadditive function such that $\Psi\leq\Phi$ holds on $[0,\ell(I)[\,$. To show that $\Psi\leq\Phi^\alpha$, let $u\in[0,\ell(I)[\,$ and $\varepsilon>0$ be arbitrary. Then 
there exist $n \in\N$ and $u_1,\dots,u_n\in\,]-\ell(I),\ell(I)[\,$ such that 
\Eq{u}
{u=u_1+\dots+u_n \qquad\text{and}\qquad
\Phi(|u_1|)+\dots+\Phi(|u_n|)<\Phi^\alpha(u)+\varepsilon. 
}
Then, due to the absolute subadditivity of $\Psi$ and \lem{2-}, we get
\Eq{*}{
  \Psi(u)\leq \Psi(|u_1|)+\dots+\Psi(|u_n|)
  \leq \Phi(|u_1|)+\dots+\Phi(|u_n|)<\Phi^\alpha(u)+\varepsilon.
}
By the arbitrariness of $\varepsilon>0$, the inequality $\Psi(u)\leq\Phi^\alpha(u)$ follows for all $u\in[0,\ell(I)[\,$, which was to be proved.

The function $\Phi^\alpha$ being subadditive, the \prp{2} implies that
$\Phi^\alpha\leq\Phi^\sigma$ also holds.
\end{proof}

The following corollaries demonstrate cases when the subadditive and the absolutely subadditive envelopes of an error function are the identically zero functions.

\Cor{A}{Let $\Phi\in\E(I)$ such that, for all $0\leq u<\ell(I)$,
\Eq{lim}{
  \inf_{n\in\N} n\Phi\Big(\frac{u}{n}\Big)=0.
} 
Then $\Phi^\sigma=\Phi^\alpha\equiv0$. In particular, for $p>1$, $\Phi_p^\sigma=\Phi_p^\alpha=0$.}

\begin{proof} If $n\in\N$ and $0\leq u<\ell(I)$, then, by the construction of $\Phi^\sigma$, 
\Eq{*}{
  \Phi^\sigma(u)\leq n\Phi\Big(\frac{u}{n}\Big).
}
Taking the infimum of the right hand side of this inequality for $n\in\N$, we get that $\Phi^\alpha(u)\leq\Phi^\sigma(u)\leq0$, which yields $\Phi^\alpha(u)=\Phi^\sigma(u)=0$.  

For the particular case, let $\Phi=\Phi_p$ for some $p>1$. Then $1-p<0$, thus
\Eq{*}{
  \inf_{n\in\N}n\Phi_p\Big(\frac{u}{n}\Big)=\inf_{n\in\N}n^{1-p}u^p
  =\lim_{n\to\infty}n^{1-p}u^p=0,
}
which shows that $\Phi_p^\sigma=\Phi_p^\alpha=0$.
\end{proof}

\Cor{B}{Let $I$ be an unbounded interval and $\Phi\in\E(I)$ such that
\Eq{lim+}{
  \lim_{v\to\infty} \Phi(v)=0.
} 
Then $\Phi^\alpha\equiv0$. In particular, for $p<0$, $\Phi_p^\alpha=0$.}

\begin{proof} Let $0<u$. Then, by the construction of $\Phi^\alpha$, for all $v>0$, we have the inequality
\Eq{*}{
  \Phi^\alpha(u)\leq \Phi(|u+v|)+\Phi(|-v|)=\Phi(u+v)+\Phi(v).
}
Upon taking the limit $v\to\infty$, the equality \eq{lim+} yields that $\Phi^\alpha(u)=0$.
\end{proof}

The next result shows that for the notions of $\Phi$-monotonicity and $\Phi$-Hölder property, the error function $\Phi$ can always be replaced by its subadditive and absolutely subadditive envelope, respectively.

\Thm{A}{Let $\Phi\in\E(I)$. Then
\Eq{*}{
  \M_\Phi(I)=\M_{\Phi^\sigma}(I) \qquad\mbox{and}\qquad 
  \H_\Phi(I)=\H_{\Phi^\sigma}(I).
}
If, in addition, $I$ is an unbounded interval, then
\Eq{*}{
  \H_\Phi(I)=\H_{\Phi^\alpha}(I).
}}

\begin{proof}
The inclusion $\M_{\Phi^\sigma}(I)\subseteq\M_{\Phi}(I)$ is a trivial consequence of the inequality $\Phi^\sigma\leq\Phi$. To prove the reversed inclusion, let $f\in\M_{\Phi}(I)$. To show that $f$ is also $\Phi^\sigma$-monotone, let $x<y$ be arbitrary elements of $I$ and $\varepsilon>0$ be arbitrary.

For $u:=y-x<\ell(I)$, there exist $n \in\N$ and $u_1,\dots,u_n\in\R_+$ such that \eq{v} holds. For the sake of convenience, let $u_0:=0$ and 
\Eq{*}{
   x_i:=x+u_0+\dots+u_i \qquad(i\in\{0,\dots,n\}).
}
Obviously, $x=x_0\leq x_1\leq \dots\leq x_n=y$. Applying the $\Phi$-monotonicity of $f$, we get that
\Eq{*}{
  f(x_{i-1})\leq f(x_{i})+\Phi(x_{i}-x_{i-1})
    =f(x_{i})+\Phi(u_i)\qquad(i\in\{1,\dots,n\}).
}
Adding up the above inequalities for $i\in\{1,\dots,n\}$ side by side, we obtain that
\Eq{*}{
  f(x)=f(x_0)\leq f(x_n)+\Phi(u_1)+\dots+\Phi(u_n)
  <f(x_n)+\Phi^\sigma(u)+\varepsilon=f(y)+\Phi^\sigma(y-x)+\varepsilon.
}
Upon taking the limit $\varepsilon\to0$, it follows that
\Eq{*}{
  f(x)\leq f(y)+\Phi^\sigma(y-x),
}
which completes the proof of the $\Phi^\sigma$-monotonicity of $f$. 

Using \prp{H0} and the above verified equality, we get
\Eq{*}{
  \H_\Phi(I)=\M_\Phi(I)\cap(-\M_\Phi(I))
  =\M_{\Phi(I)^\sigma}\cap(-\M_{\Phi(I)^\sigma})=\H_{\Phi^\sigma}(I).
}

Now assume that $I$ is unbounded.
The inclusion $\H_{\Phi^\alpha}(I)\subseteq\H_{\Phi}(I)$ is a trivial consequence of the inequality $\Phi^\alpha\leq\Phi$. To prove the reversed inclusion, let $f\in\H_{\Phi}(I)$. To show that $f$ is also $\Phi^\alpha$-H\"older, let $x,y\in I$ and $\varepsilon>0$ be arbitrary. We may assume that $x<y$. Then $u:=y-x<\ell(I)$ and there exist $n\in\N$ and $u_1,\dots,u_n\in\,]-\ell(I),\ell(I)[\,$ such that \eq{u} holds. 

Now we have to distinguish two cases according to the unboundedness of $I$.
Assume first that $I$ is unbounded from below. Then we may assume that $-\ell(I):=u_0<u_1\leq \dots\leq u_n<\ell(I)$. In view of \eq{u}, we have that $u_n>0$, therefore there exists a unique $k\in\{1,\dots,n\}$ such that $u_{k-1}\leq0<u_{k}$. For the sake of convenience, let 
\Eq{*}{
  x_0:=x,\qquad x_i:=x+u_1+\dots +u_i\qquad(i\in\{1,\dots,n\}).
}
By the construction of $k$, it follows that
\Eq{*}{
   x=x_0\geq x_1\geq\dots\geq x_{k-1}<x_k\leq\dots\leq x_n=y.
}
Therefore, $x_i\leq\max(x,y)=$ for all $i\in\{1,\dots,n\}$. Thus the unboundedness of $I$ from below yields that $x_1,\dots,x_n\in I$ hold. Applying the $\Phi$-Hölder property of $f$, we get that
\Eq{*}{
  f(x_{i-1})\leq f(x_{i})+\Phi(|x_{i}-x_{i-1}|)
    =f(x_{i})+\Phi(|u_i|)\qquad(i\in\{1,\dots,n\}).
}
Adding up the above inequalities for $i\in\{1,\dots,n\}$ side by side, we obtain that
\Eq{*}{
  f(x)=f(x_0)\leq f(x_n)+\Phi(|u_1|)+\dots+\Phi(|u_n|)
  <f(x_n)+\Phi^\alpha(u)+\varepsilon=f(y)+\Phi^\alpha(y-x)+\varepsilon.
}
Upon taking the limit $\varepsilon\to0$, it follows that
\Eq{xy}{
  f(x)\leq f(y)+\Phi^\alpha(|y-x|).
}
In the case when $I$ is unbounded from above, one should take the ordering $-\ell(I)<u_n\leq \dots\leq u_1<u_0:=\ell(I)$ and use a completely similar argument to obtain inequality \eq{xy}.

Finally, interchanging the roles of $x$ and $y$ in the above proof, we can get 
\Eq{*}{
  f(y)\leq f(x)+\Phi^\alpha(|y-x|),
} 
which, together with \eq{xy}, shows that $f$ is $\Phi^\alpha$-H\"older and hence completes the proof.
\end{proof}

\Cor{AA}{Let $\Phi\in\E(I)$ such that, for all $0\leq u<\ell(I)$, \eq{lim} holds. Then $\M_\Phi(I)$ equals the class of increasing functions on $I$ and $\H_\Phi(I)$ consists of constant functions.}

\begin{proof} In view of \cor{A}, we have that $\Phi^\sigma=\Phi^\alpha\equiv0$. Combining this with the result of the \thm{A}, we get that $\M_\Phi(I)=\M_0(I)$ and $\H_\Phi(I)=\H_0(I)$, which is equivalent to the statement.
\end{proof}

The following result demonstrates that the subadditive and increasing error functions optimally determine the corresponding classes of monotone and Hölder functions.

\Thm{Opt}{Let $\Phi\in\E(I)$ be a increasing and subadditive function with $\Phi(0)=0$ and let $\Psi\in\E(I)$ satisfy the inequality $\Psi\leq\Phi$. Then $\M_\Psi(I)=\M_\Phi(I)$ if and only if $\Psi=\Phi$. Similarly, $\H_\Psi(I)=\H_\Phi(I)$ if and only if $\Psi=\Phi$.}

\begin{proof} The inclusions $\M_\Psi(I)\subseteq\M_\Phi(I)$ and $\H_\Psi(I)\subseteq\H_\Phi(I)$ follow from the inequality $\Psi\leq\Phi$. To prove the statement, it suffices to show that if $\Psi(p)<\Phi(p)$ for some $p\in\,]0,\ell(I)[\,$ then both inclusions $\M_\Psi(I)\subseteq\M_\Phi(I)$ and $\H_\Psi(I)\subseteq\H_\Phi(I)$ are proper. For this, we construct a function $f:I\to\R$ such that $f\in\H_\Phi(I)$ but $f\notin\M_\Psi(I)$.

The inclusion $p\in\,]0,\ell(I)[\,$ implies that there exist $u,v\in I$ such that $p=v-u$. Define 
\Eq{*}{
   f(x):=\begin{cases}
         0&\mbox{if } x\leq u\\
         -\Phi(x-u)&\mbox{if } u<x
         \end{cases}\qquad (x\in I).
}
Then $-f$ is increasing because $\Phi$ is increasing, hence $-f\in\M_\Phi(I)$. To prove that $f\in\M_\Phi(I)$, we fix $x,y\in I$ with $x\leq y$ and distinguish three cases. 

If $x\leq y\leq u$, then $f(x)=f(y)=0$, hence the inequality \eq{H1} is a consequence of the nonnegativity of $\Phi$. 

If $x\leq u<y$, then $f(x)=0$ and $f(y)=-\Phi(y-u)$, therefore the inequality 
\eq{H1} is now equivalent to
\Eq{*}{
  \Phi(y-u)\leq \Phi(y-x),
}
which is a consequence of the increasingness of $\Phi$.

If $u<x\leq y$, then $f(x)=-\Phi(x-u)$ and $f(y)=-\Phi(y-u)$, therefore the inequality \eq{H1} is now equivalent to
\Eq{*}{
  \Phi(y-u)\leq \Phi(y-x)+\Phi(x-u),
}
which is a consequence of the subadditivity of $\Phi$.

This completes the proof of the inclusion $f\in\M_\Phi(I)$ and hence shows that $f\in\H_\Phi(I)$. To complete the proof, we have to verify that $f\notin \M_\Psi(I)$. Indeed, we have that
\Eq{*}{
  f(u)-f(v)=\Phi(v-u)=\Phi(p)>\Psi(p)=\Psi(v-u).
}
This strict inequality shows that $f$ cannot be $\Psi$-monotone.
\end{proof}

\section{$\Phi$-monotone and $\Phi$-Hölder envelopes}

As we have seen it in \prp{MC2} and \prp{H2}, the classes $\M_\Phi(I)$ and $\H_\Phi(I)$ are closed with respect to pointwise infimum and maximum. Therefore, for any function $f:I\to\R$, the supremum of all $\Phi$-monotone ($\Phi$-Hölder) functions below $f$ (provided that there is at least one such function) is the largest $\Phi$-monotone ($\Phi$-Hölder) function which is smaller than or equal to $f$. Similarly, the infimum of all $\Phi$-monotone ($\Phi$-Hölder) functions above $f$ (provided that there is at least one such function) is the smallest $\Phi$-monotone ($\Phi$-Hölder) function which is bigger than or equal to $f$. The next result offers a formula for these enveloping functions.

\Prp{nv}{Let $\Phi\in\E(I)$ with $\Phi(0)=0$ and let $f:I\to \R$ be a function which admits a $\Phi$-monotone minorant. Then the function $\underline{M}_\Phi(f)$ defined by
\Eq{*}{
  \underline{M}_\Phi(f)(x):=\inf_{x\leq y} \big(f(y)+\Phi^\sigma(y-x)\big) \qquad(x\in I)
}
is real-valued and is the largest $\Phi$-monotone function which is smaller than or equal to $f$. Analogously, if $f$ admits a $\Phi$-monotone majorant, then the function $\overline{M}_\Phi(f)$ defined by
\Eq{*}{
  \overline{M}_\Phi(f)(x):=\sup_{y\leq x} \big(f(y)-\Phi^\sigma(x-y)\big) \qquad(x\in I)
}
is real-valued and is the smallest $\Phi$-monotone function which is bigger than or equal to $f$.}

\begin{proof}
Obviously, $\underline{M}_\Phi(f)$ cannot take the value $+\infty$ at any point in $I$, i.e., $\underline{M}_\Phi(f)(x)<+\infty$ for all $x\in I$. The condition $\Phi(0)=0$ implies that $\Phi^\sigma(0)=0$, therefore, by taking $y=x$ in the defining  formula of $\underline{M}_\Phi(f)(x)$, we get that $\underline{M}_\Phi(f)(x)\leq f(x)$ holds for all $x\in I$. 

Now suppose $g$ is a $\Phi$-monotone function such that $g\leq f$ holds (by the assumption, there is at least one such function $g$). Then, by \thm{A}, $g$ is also $\Phi^\sigma$-monotone. In order to show that $g\leq \underline{M}_\Phi(f)$, let $x\in I$ be arbitrarily fixed. Then, for all $y\in I$ with $x\leq y$, we have
\Eq{*}{
 g(x)\leq g(y)+\Phi^\sigma(y-x)
 \leq f(y)+\Phi^\sigma(y-x).
}
Upon taking the infimum of the right hand side with respect to $y\geq x$, we get
\Eq{*}{
 g(x) \leq \inf_{x\leq y}\big(f(y)+\Phi^\sigma(y-x)\big)=\underline{M}_\Phi(f)(x),
}
which proves the desired inequality $g(x)\leq \underline{M}_\Phi(f)(x)$ and also that $\underline{M}_\Phi(f)$ cannot take the value $-\infty$ at any point of $I$.

To see that $\underline{M}_\Phi(f)$ itself is $\Phi$-monotone, it is sufficient to show that $\underline{M}_\Phi(f)$ is $\Phi^\sigma$-monotone. Let $u,v\in I$ with $u\leq v$. Then, using the subadditivity of $\Phi^\sigma$, we obtain
\Eq{*}{
  \underline{M}_\Phi(f)(u)&=\inf_{u\leq y} \big(f(y)+\Phi^\sigma(y-u)\big)
  \leq \inf_{v\leq y} \big(f(y)+\Phi^\sigma(y-u)\big) \\
  &\leq \inf_{v\leq y} \big(f(y)+\Phi^\sigma(y-v)\big)+\Phi^\sigma(v-u)
  = \underline{M}_\Phi(f)(v)+\Phi^\sigma(v-u),
}
which completes the proof of the $\Phi$-monotonicity of $\underline{M}_\Phi(f)$.

The proof of the second assertion is completely similar.
\end{proof}

The following result is of a sandwich type one.

\Cor{MSW}{Let $\Phi\in\E(I)$ with $\Phi(0)=0$ and let $g,h:I\to\R$. Then in order that there exist a $\Phi$-monotone function $f:I\to\R$ between $g$ and $h$ it is necessary and sufficient that, for all $x,y\in I$ with $x\leq y$, the inequality
\Eq{gh}{
  g(x)\leq h(y)+\Phi^\sigma(y-x)
}
be valid.}

\begin{proof}
Assume first that $f$ is a $\Phi$-monotone function such that $g\leq f\leq h$.
Then, $f$ is $\Phi^\sigma$-monotone and, for all $x,y\in I$ with $x\leq y$, we have
\Eq{*}{
  g(x)\leq f(x)\leq f(y)+\Phi^\sigma(y-x)\leq h(y)+\Phi^\sigma(y-x),
}
i.e., \eq{gh} holds.

Conversely, assume that \eq{gh} holds true for all $x,y\in I$ with $x\leq y$.
For a fixed $x\in I$, define
\Eq{*}{
  f(x):=\underline{M}_\Phi(h)(x)=\inf_{x\leq y} \big(h(y)+\Phi^\sigma(y-x)\big).
}
Now, in view of inequality \eq{gh}, we have that $g(x)\leq f(x)$. By taking $y=x$ in the definition of $f$, the condition $\Phi(0)=0$ ensures that $f(x)\leq h(x)$ is also valid. Finally, arguing similarly as at the end of the proof of \prp{nv},
it follows that $f$ is $\Phi^\sigma$-monotone and hence $\Phi$-monotone as well.
\end{proof}

\Prp{H-env}{Let $I$ be an unbounded interval, $\Phi\in\E(I)$ with $\Phi(0)=0$  and let $f:I\to \R$ be a function which admits a $\Phi$-Hölder minorant. Then the function $\underline{H}_\Phi(f)$ defined by
\Eq{*}{
  \underline{H}_\Phi(f)(x):=\inf_{y\in I} \big(f(y)+\Phi^\alpha(|y-x|)\big) \qquad(x\in I)
}
is real-valued and is the largest $\Phi$-Hölder function which is smaller than or equal to $f$. Analogously, if $f$ admits a $\Phi$-Hölder majorant, then the function $\overline{H}_\Phi(f)$ defined by
\Eq{*}{
  \overline{H}_\Phi(f)(x):=\sup_{y\in I} \big(f(y)-\Phi^\alpha(|x-y|)\big) \qquad(x\in I)
}
is real-valued and is the smallest $\Phi$-Hölder function which is bigger than or equal to $f$.}

\begin{proof}
Obviously, $\underline{H}_\Phi(f)$ cannot take the value $+\infty$ at any point in $I$, i.e., $\underline{H}_\Phi(f)(x)<\infty$ for all $x\in I$. The condition $\Phi(0)=0$ implies that $\Phi^\alpha(0)=0$, therefore, by taking $y=x$ in the defining  formula of $\underline{H}_\Phi(f)(x)$, we get that $\underline{H}_\Phi(f)(x)\leq f(x)$ holds for all $x\in I$. 

Now suppose $g$ is a $\Phi$-H\"older function such that $g\leq f$ holds (by the assumption, there is at least one such function $g$). Then, by \thm{A}, $g$ is also $\Phi^\alpha$-H\"older. In order to show that $g\leq \underline{H}_\Phi(f)$, let $x\in I$ be arbitrarily fixed. Then, for all $y\in I$, we have
\Eq{*}{
 g(x)\leq g(y)+\Phi^\alpha(|y-x|)
 \leq f(y)+\Phi^\alpha(|y-x|),
}
Upon taking the infimum of the right hand side with respect to $y\in I$, we get
\Eq{*}{
 g(x) \leq \inf_{ y}\big(f(y)+\Phi^\alpha(|y-x|)\big)=\underline{H}_\Phi(f)(x),
}
which proves the desired inequality $g(x)\leq \underline{H}_\Phi(f)(x)$ and also that $\underline{H}_\Phi(f)$ cannot take the value $-\infty$ at any point of $I$.

To see that $\underline{H}_\Phi(f)$ itself is $\Phi$-H\"older, it is sufficient to show that $\underline{H}_\Phi(f)$ is $\Phi^\alpha$-H\"older. For any  $u,v\in I$, using the absolute subadditivity of $\Phi^\alpha$, we obtain
\Eq{*}{
  \underline{H}_\Phi(f)(u)&=\inf_{ y} \big(f(y)+\Phi^\alpha(|y-u|\big)
  \leq \inf_{ y} \bigg(f(y)+\Phi^\alpha\big(|(y-v)|\big)+\Phi^\alpha\big(|(v-u)|\big)\bigg) \\
  &= \inf_{ y} \Big(f(y)+\Phi^\alpha(|y-v|)\Big)+\Phi^\alpha(|v-u|)
  = \underline{H}_\Phi(f)(v)+\Phi^\alpha(|v-u|).
}
In the same pattern as above, interchanging the roles of $u$  and $v$ in the above equation, we will obtain
\Eq{*}{
  \underline{H}_\Phi(f)(v)\leq\underline{H}_\Phi(f)(u)+\Phi^\alpha(|v-u|),
}
which shows that  $\underline{H}_\Phi(f)$ is $\Phi$-H\"older.

The proof of the second assertion is completely similar.
\end{proof}

\Cor{HSW}{Let $I$ be an unbounded interval, let $\Phi\in\E(I)$ with $\Phi(0)=0$ and let $g,h:I\to\R$. Then in order that there exist a $\Phi$-Hölder function $f:I\to\R$ between $g$ and $h$ it is necessary and sufficient that, for all $x,y\in I$, the inequality
\Eq{gh+}{
  g(x)\leq h(y)+\Phi^\alpha(|y-x|)
}
be valid.}

\begin{proof}
Assume first that $f$ is a $\Phi$-Hölder function such that $g\leq f\leq h$.
Then, $f$ is $\Phi^\alpha$-Hölder and, for all $x,y\in I$, we have
\Eq{*}{
  g(x)\leq f(x)\leq f(y)+\Phi^\alpha(|y-x|)\leq h(y)+\Phi^\alpha(|y-x|),
}
i.e., \eq{gh+} holds.

Conversely, assume that \eq{gh+} holds true for all $x,y\in I$.
For a fixed $x\in I$, define
\Eq{*}{
  f(x):=\underline{H}_\Phi(h)(x)=\inf_{y\in I} \big(h(y)+\Phi^\alpha(|y-x|)\big).
}
Now, in view of inequality \eq{gh+}, we have that $g(x)\leq f(x)$. By taking $y=x$ in the definition of $f$, the conditions $\Phi(0)=0$ ensures that $f(x)\leq h(x)$ is also valid. Finally, arguing similarly as at the end of the proof of \prp{H-env}, it follows that $f$ is $\Phi^\alpha$-Hölder and hence $\Phi$-Hölder as well.
\end{proof}

Before we formulate and prove the next theorem we shall need the following auxiliary result.

\Lem{PP}{Let $\Phi,\Psi\in\E(I)$ such that $(-\Phi)$ is $\Psi$-monotone on $]0,\ell(I)[$. Then $(-\Phi^\sigma)$ is also $\Psi$-monotone on $]0,\ell(I)[$.}

\begin{proof}
To prove this lemma, let $x,y\in ]0,\ell(I)[$  with $x<y$ and $\varepsilon>0$ be arbitrary. Then by definition of $\Phi^\sigma$, there exist $n\in N$ and $u_1,\dots,u_n\in\R_+$ such that $x=u_1+\dots +u_n$ satisfying
\Eq{mmm}
{\Phi(u_1)+\dots+\Phi(u_n)<\Phi^\sigma(x)+\varepsilon.
}
Using the $\Psi$-monotonicity of $(-\Phi)$, we have
\Eq{*}
{(-\Phi)(u_n)\leq (-\Phi)(u_n+(y-x))+\Psi(y-x),}
from which we obtain 
\Eq{*}
{\Phi(u_n+(y-x))\leq \Phi(u_n)+\Psi(y-x).
}
Observe that $y=u_1+\dots+u_{n-1}+(u_n+(y-x))$. Thus, using the inequality in \eq{mmm}, we arrive at
\Eq{*}
{\Phi^\sigma(y)&\leq\Phi(u_1)+\dots+\Phi(u_{n-1})+\Phi(u_n+(y-x))\\
&\leq\Phi(u_1)+\dots+\Phi(u_{n-1})+\Phi(u_n)+\Psi(y-x)
<\Phi^\sigma(x)+\Psi(y-x)+\varepsilon.
}
As $\varepsilon$ is an arbitrary positive number, we can conclude that $(-\Phi^\sigma)(x)\leq (-\Phi^\sigma)(y)+\Psi(y-x)$, which completes the proof of the $\Psi$ monotonicity of $(-\Phi^\sigma)$.
\end{proof}

\Thm{PP}{Let $\Phi,\Psi\in\E(I)$ such that $(-\Phi)$ is $\Psi$-monotone on $]0,\ell(I)[$ and let $f:I\to\R$. Then $f$ is $\Phi$-monotone if and only if there exist two $\Psi$-monotone functions $f_*,f^*:I\to \R$ such that $f_*\leq f\leq f^*$ hold on $I$ and, for all $x,y\in I$ with  $x<y$,
\Eq{gf}{
f(x)\leq f_*(y)+\Phi^\sigma(y-x) \qquad\mbox{and}\qquad 
f^*(x)\leq f(y)+\Phi^\sigma(y-x).
}}

\begin{proof}
First assume that $f$ is $\Phi$-monotone. Then, by \thm{A}, it is also $\Phi^\sigma$-monotone. For a fixed point $x\in I$, define
\Eq{*}{
f_*(x):=\sup_{u<x}\Big(f(u)-\Phi^\sigma(x-u)\Big)\qquad \mbox{and} \qquad 
f^*(x):=\inf_{x<v}\Big(f(v)+\Phi^\sigma(v-x)\Big).
}
In view of the $\Phi^\sigma$-monotonicity of $f$, for all $u<x<v$, we have that
\Eq{*}{
  f(u)-\Phi^\sigma(x-u)\leq f(x) \qquad\mbox{and}\qquad 
  f(x)\leq f(v)+\Phi^\sigma(v-x).
}
Therefore, upon taking the supremum for $u<x$ and the infimum for $x<v$, we get
that $f_*(x)\leq f(x)$ and $f(x)\leq f^*(x)$, respectively. That is, we have that $f_*$ and $f^*$ are real valued functions and the inequalities $f_*\leq f\leq f^*$ hold on $I$.

In the next step, we establish the $\Psi$-monotonicity of $f_*$ and $f^*$. 
By the definition of $f^*$, for all $x,y\in I$ with $x<y$, we have
\Eq{g**}{
  f_*(x)=\sup_{u<x}\Big(f(u)-\Phi^\sigma(x-u)\Big)
        \leq \sup_{u<y}\Big(f(u)-\Phi^\sigma(x-u)\Big).
}
By \lem{PP}, we have the $\Psi$-monotonicity of $(-\Phi^\sigma)$, which, for all $u\in I$ with $u<x$, implies
\Eq{*}{
  -\Phi^\sigma(x-u)\leq -\Phi^\sigma(y-u)+\Psi((y-u)-(x-u))
  =-\Phi^\sigma(y-u)+\Psi(y-x).
}
Applying this inequality to the right most expression of inequality \eq{g**}, we arrive at
\Eq{*}{
  f_*(x)\leq \sup_{u<y}\Big(f(u)-\Phi^\sigma(x-u)\Big)
  \leq \sup_{u<y}\Big(f(u)-\Phi^\sigma(y-u)\Big)+\Psi(y-x)=f_*(y)+\Psi(y-x),
}
which shows that $f_*$ is also $\Psi$-monotone.

Take $x,y\in I$ with $x<y$. Then, by the definition of $f^*$, we have
\Eq{g*}{
  f^*(x)=\inf_{x<v}\Big(f(v)+\Phi^\sigma(v-x)\Big)
        \leq \inf_{y<v}\Big(f(v)+\Phi^\sigma(v-x)\Big).
}
By the $\Psi$-monotonicity of $(-\Phi^\sigma)$, for all $v\in I$ with $y<v$, we obtain
\Eq{*}{
  \Phi^\sigma(v-x)\leq\Phi^\sigma(v-y)+\Psi((v-x)-(v-y))
  =\Phi^\sigma(v-y)+\Psi(y-x).
}
Applying this inequality to the right most expression of inequality \eq{g*}, we arrive at
\Eq{*}{
  f^*(x)\leq \inf_{y<v}\Big(f(v)+\Phi^\sigma(v-x)\Big)
  \leq \inf_{y<v}\Big(f(v)+\Phi^\sigma(v-y)\Big)+\Psi(y-x)=f^*(y)+\Psi(y-x),
}
which proves that $f^*$ is $\Psi$-monotone.

Finally, for $x,y\in I$ with $x<y$, from the definitions of $f_*$ and $f^*$, we obtain the inequalities 
\Eq{*}{
  f(x)-\Phi^\sigma(y-x)\leq f_*(y) \qquad\mbox{and}\qquad 
  f^*(x)\leq f(y)+\Phi^\sigma(y-x),
}
respectively, which prove that $f_*$ and $f^*$ satisfy the inequalities stated in \eq{gf}.

Conversely, if the first inequality in \eq{gf} holds for some function $f_*:I\to\R$ satisfying $f_*\leq f$, then $f(x)\leq f_*(y)+\Phi^\sigma(y-x)\leq f(y)+\Phi^\sigma(y-x)$, which shows the $\Phi^\sigma$-monotonicity of $f$. Similarly, the existence of a function $f^*:I\to\R$ satisfying $f\leq f^*$ and the second inequality of \eq{gf}, also implies that $f$ is $\Phi^\sigma$-monotone.
\end{proof}

By taking the error function $\Psi\equiv0$, the previous theorem directly implies the following result. Observe that, in this case, $\Psi$-monotonicity is equivalent to increasingness.

\Cor{PP}{Let $\Phi\in\E(I)$ such that $\Phi$ is decreasing on $]0,\ell(I)[$ and let $f:I\to\R$. Then, $f$ is $\Phi$-monotone if and only if there exist two increasing functions $f_*,f^*:I\to \R$ such that $f_*\leq f\leq f^*$ hold on $I$ and, for all $x,y\in I$ with $x<y$, the inequalities in \eq{gf} are satisfied.}

The analogue of \lem{PP} for the $\Psi$-Hölder setting is contained in the following lemma.

\Lem{qq}{Let $I$ be an unbounded interval, let $\Phi,\Psi\in\E(I)$ such that $\Phi\circ|\cdot|$ is $\Psi$-Hölder on $\R$. Then $\Phi^\alpha\circ|\cdot|$ is also $\Psi$-Hölder on $\R$.}

\begin{proof}
To prove this lemma, let $x,y\in\R$ and $\varepsilon>0$ be arbitrary. Then by definition of $\Phi^\alpha$, there exist $n\in N$ and $u_1,\dots,u_n\in\R$ with $y=u_1+\dots +u_n$ satisfying
\Eq{mm}
{\Phi(|u_1|)+\dots+\Phi(|u_n|)<\Phi^\alpha(|y|)+\varepsilon.
}
Using the $\Psi$-Hölder property of $\Phi\circ|\cdot|$, we have
\Eq{*}
{\Phi(|u_n+(x-y)|)\leq \Phi(|u_n|)+\Psi(|y-x|).
}
Observe that $x=u_1+\dots+u_{n-1}+(u_n+(x-y))$. Thus, using the inequality in \eq{mm}, we arrive at
\Eq{*}
{\Phi^\alpha(|x|)&\leq\Phi(|u_1|)+\dots+\Phi(|u_{n-1}|)+\Phi(|u_n+(x-y)|)\\
&\leq\Phi(|u_1|)+\dots+\Phi(|u_{n-1}|)+\Phi(|u_n|)+\Psi(|y-x|)
<\Phi^\alpha(|y|)+\Psi(|y-x|)+\varepsilon.
}
As $\varepsilon$ is an arbitrary positive number, we can conclude that $\Phi^\alpha\circ|\cdot|$ is $\Psi$-Hölder.
\end{proof}

\Rem{qq}{The property that $\Phi\circ|\cdot|$ is $\Psi$-Hölder on $\R$ means that, for all $x,y\in\R$,
\Eq{*}{
  \Phi(|x|)\leq\Phi(|y|)+\Psi(|y-x|).
}
Denoting $|x|$ and $|y|$ by $u$ and $v$, respectively, the above inequality implies, for all $u,v\geq0$,
\Eq{*}{
  \Phi(u)\leq\Phi(v)+\min(\Psi(|v-u|),\Psi(u+v)).
}
Conversely, one can see that the $\Psi$-Hölder property of the function $\Phi\circ|\cdot|$ is a consequence of the last inequality. Therefore, if $\Psi$ is increasing, then $\Phi\circ|\cdot|$ is $\Psi$-Hölder if and only if $\Phi$ is $\Psi$-Hölder.}

\Thm{qq}{Let $I$ be an unbounded interval, $\Phi,\Psi\in\E(I)$ such that $\Phi\circ|\cdot|$ is $\Psi$-H\"older on $\R$. Then $f:I\to\R$ is $\Phi$-H\"older if and only if there exist two $\Psi$-H\"older functions $f_*,f^*:I\to \R$ such that $f_*\leq f\leq f^*$ holds on $I$ and, for all $x,y\in I $ with $x\neq y$
\Eq{gf+}{
f(x)\leq f_*(y)+\Phi^\alpha(|y-x|) \qquad\mbox{and}\qquad 
f^*(x)\leq f(y)+\Phi^\alpha(|y-x|).
}
Additionally, for all $x\in I$,
\Eq{dif}
{f^*(x)-f_*(x)\leq\inf_{y\in I}2\Phi(|y-x|).
}}

\begin{proof}
First assume that $f$ is $\Phi$-H\"older. Then, by \thm{A}, $f$ is also $\Phi^\alpha$-H\"older. For a fixed point $x\in I$, define
\Eq{*}{
f_*(x):=\sup_{u\in I}\Big(f(u)-\Phi^\alpha(|u-x|)\Big)\qquad \mbox{and} \qquad 
f^*(x):=\inf_{u\in I}\Big(f(u)+\Phi^\alpha(|u-x|)\Big).
}
In view of the $\Phi^\alpha$-H\"older property of $f$, for all $u,x\in I$, we have that
\Eq{*}{
  f(u)-\Phi^\alpha(|u-x|)\leq f(x) \qquad\mbox{and}\qquad f(x)\leq f(u)+\Phi^\alpha(|u-x|).
}
Therefore, upon taking the supremum and infimum for all $u\in I$, we get that $f_*(x)\leq f(x)$ and $f(x)\leq f^*(x)$, respectively. That is, $f_*$ and $f^*$ are real valued functions and $f_*\leq f\leq f^*$ is satisfied on $I$.
In the next step, we establish that $f_*$ and $f^*$ are $\Psi$-H\"older. 

In view of \lem{qq}, the $\Psi$-H\"older property of $\Phi\circ|\cdot|$ yields that $\Phi^\alpha\circ|\cdot|$ is also $\Psi$-Hölder. Therefore, for all $u,x,y \in I$, we have
\Eq{*}{
  -\Phi^\alpha(|x-u|)\leq -\Phi^\alpha(|y-u|)+\Psi(|(y-u)-(x-u)|)\leq-\Phi^\alpha(|u-y|)+\Psi(|y-x|).
}
Applying this inequality to the definition of $f_*$, we arrive at
\Eq{*}{
  f_*(x)=\sup_{u\in I} \Big(f(u)-\Phi^\alpha(|u-x|)\Big)
  \leq \sup_{u\in I}\Big(f(u)-\Phi^\alpha(|u-y|)\Big)+\Psi(|y-x|)=f_*(y)+\Psi(|y-x|),
}
which shows that $f_*$ is $\Psi$-H\"older.

Similarly, by the $\Psi$-H\"older property of $\Phi^\alpha\circ|\cdot|$, for all $u,x,y\in I$, we obtain
\Eq{*}{
  \Phi^\alpha(|u-x|)\leq \Phi^\alpha(|u-y|)+\Psi(\left|(u-y)-(u-x)\right|)=\Phi^\alpha(|u-y|)+\Psi(|y-x|).
}
Applying this inequality to the definition of $f^*$, it follows that
\Eq{*}{
  f^*(x)=\inf_{u\in I}\Big(f(u)+\Phi^\alpha(|u-x|)\Big)
  \leq \inf_{u\in I}\Big(f(u)+\Phi^\alpha(|u-y|)\Big)+\Psi(|y-x|)=f^*(y)+\Psi(|y-x|),
}
which proves that $f^*$ is also $\Psi$-H\"older.

Next we prove that $f_*$ and $f^*$ satisfy the inequalities stated in \eq{gf+}. Indeed, for $x,y\in I$, from the definitions of $f_*$ and $f^*$, we obtain the inequalities 
\Eq{*}{
  f(x)-\Phi^\alpha(|y-x|)\leq f_*(y) \qquad\mbox{and}\qquad 
  f^*(x)\leq f(y)+\Phi^\alpha(|y-x|),
}
respectively. Conversely, if the first inequality in \eq{gf+} holds for some function $f_*:I\to\R$ satisfying $f_*\leq f$, then $f(x)\leq f_*(y)+\Phi^\alpha(|y-x|)\leq f(y)+\Phi^\alpha(|y-x|)$, which shows that $f$ is $\Phi^\alpha$-H\"older. Similarly, the existence of a function $f^*:I\to\R$ satisfying $f\leq f^*$ and the second inequality of \eq{gf+}, also implies that $f$ is $\Phi^\alpha$-H\"older.

Finally, to obtain the last inequality \eq{dif} of \thm{qq}, we interchange $x$ and $y$ in the first inequality of \eq{gf+} and we obtain 
\Eq{*}
{-f_*(x)\leq -f(y)+\Phi(|y-x|)
}
By summing up this inequality with the second inequality of \eq{gf+} side by side, we will reach at our desired conclusion.
\end{proof}

\section{Jordan-type decomposition of functions with bounded $\Phi$-variation}

Let $\Phi\in\E(I)$. Then a function $f:I\to\R$ is called \emph{delta-$\Phi$-monotone} if it is the difference of two $\Phi$-monotone functions. In what follows, we shall extend the celebrated Jordan Decomposition Theorem for delta-$\Phi$-monotone functions. For this purpose, we have to extend the notion of total variation to this more general setting.

Let $[a,b]\subseteq I$ and let $\tau=(t_0,\dots,t_n)$ be a partition of the interval $[a,b]$ (i.e., $a=t_0<t_1<\dots<t_n=b$). Then the \emph{$\Phi$-variation of $f$ with respect to $\tau$} is defined by
\Eq{*}{
  V^\Phi(f;\tau):=\sum_{i=1}^n \big(|f(t_i)-f(t_{i-1})|-\Phi(t_i-t_{i-1})\big).
}
Finally, the \emph{total $\Phi$-variation of $f$ on the interval $[a,b]$} is defined by
\Eq{*}{
 V^\Phi_{[a,b]}f:=\sup\{V^\Phi(f;\tau)\mid \tau\mbox{ is a partition of }[a,b]\}.
}

\Lem{JDT}{Let $\Phi\in\E(I)$. Then, for all $f:I\to\R$ and $a<b<c$ in $I$, we have
\Eq{JDT}{
  V^\Phi_{[a,b]}f+V^\Phi_{[b,c]}f\leq V^\Phi_{[a,c]}f.
}}

\begin{proof} Let 
\Eq{*}{
  u<V^\Phi_{[a,b]}f \qquad\mbox{and}\qquad v<V^\Phi_{[b,c]}f
}
be arbitrary real numbers. Then there exist a partition $\tau=(t_0,\dots,t_n)$ of $[a,b]$ and partition $\sigma=(s_0,\dots,s_m)$ of $[b,c]$ such that
\Eq{*}{
  u<\sum_{i=1}^n \big(|f(t_i)-f(t_{i-1})|-\Phi(t_i-t_{i-1})\big) \qquad\mbox{and}\qquad
  v<\sum_{j=1}^m \big(|f(s_j)-f(s_{j-1})|-\Phi(s_j-s_{j-1})\big).
}
Observe that $\tau\cup\sigma:=(t_0,\dots,t_n=b=s_0,\dots,s_m)$ is a partition of the interval $[a,c]$. Therefore, adding the above inequalities side by side, we get
\Eq{*}{
  u+v<V^\Phi(f;\tau\cup\sigma)\leq V^\Phi_{[a,c]}f.
}
Using the arbitrariness of $u$ and $v$, it follows that \eq{JDT} holds.
\end{proof}

Our first result characterizes those functions whose total $\Phi$-variation is nonpositive on every subinterval of $I$.

\Thm{Hold}{Let $\Phi\in\E(I)$. Then $V^\Phi_{[a,b]}f\leq0$ holds for all $a<b$ in $I$ if and only if $f$ is a $\Phi$-Hölder function.}

\begin{proof} Assume first that $f$ is a $\Phi$-Hölder function and let $a<b$ in $I$. Then, for any partition $\tau=(t_0,\dots,t_n)$ of $[a,b]$, the $\Phi$-Hölder property of $f$ yields
\Eq{*}{
  |f(t_i)-f(t_{i-1})|-\Phi(t_i-t_{i-1})\leq 0 \qquad (i\in\{1,\dots,n\}).
}
After summation, this results that $V^\Phi(f;\tau)\leq0$ for all partition $\tau$ and hence $V^\Phi_{[a,b]}f\leq0$.

Now assume that, for all $a<b$ in $I$, $V^\Phi_{[a,b]}f\leq0$. Then $V^\Phi(f;\tau)\leq0$, where $\tau$ is the trivial partition $t_0=a$, $t_1=b$. Therefore,
\Eq{*}{
  |f(b)-f(b)|-\Phi(b-a)\leq 0.
}
This shows that $f$ is $\Phi$-Hölder, indeed.
\end{proof}

The main results of this section are as follows.

\Thm{JDT1}{Let $\Phi,\Psi\in\E(I)$. If $f:I\to\R$ is the difference of a $\Phi$-monotone and a $\Psi$-monotone functions, then the total $2\max(\Phi,\Psi)$-variation of $f$ is finite on every compact subinterval of $I$.}

\begin{proof} Assume that $f=g-h$, where $g:I\to\R$ is $\Phi$-monotone and $h:I\to\R$ is $\Psi$-monotone. Let $[a,b]\subseteq I$ and let $\tau=(t_0,\dots,t_n)$ be a partition of $[a,b]$. Then, by the monotonicity properties of $g$ and $h$, for all $i\in\{1,\dots,n\}$, we have
\Eq{*}{
  g(t_i)-g(t_{i-1})+\Phi(t_i-t_{i-1})\geq0 \qquad\mbox{and}\qquad
  h(t_i)-h(t_{i-1})+\Psi(t_i-t_{i-1})\geq0.
}
Therefore, by the triangle inequality,
\Eq{*}{
  |f(t_i)-f(t_{i-1})|&-2\max(\Phi,\Psi)(t_i-t_{i-1})\\
  &=\big|[g(t_i)-g(t_{i-1})+\Phi(t_i-t_{i-1})]-[h(t_i)-h(t_{i-1})+\Psi(t_i-t_{i-1})] \\
  &\qquad  +(\Psi-\Phi)(t_i-t_{i-1})\big|-2\max(\Phi,\Psi)(t_i-t_{i-1})\\
  &\leq[g(t_i)-g(t_{i-1})+\Phi(t_i-t_{i-1})]+[h(t_i)-h(t_{i-1})+\Psi(t_i-t_{i-1})]|\\
  &\qquad +|\Psi-\Phi|(t_i-t_{i-1})-2\max(\Phi,\Psi)(t_i-t_{i-1}) \\
  &=g(t_i)-g(t_{i-1})+h(t_i)-h(t_{i-1})+(\Phi+\Psi+|\Psi-\Phi|-2\max(\Phi,\Psi))(t_i-t_{i-1})\\
  &=g(t_i)-g(t_{i-1})+h(t_i)-h(t_{i-1}).
}
Summing up these inequalities side by side for $i\in\{1,\dots,n\}$, we obtain
\Eq{*}{
  V^{2\max(\Phi,\Psi)}(f;\tau)\leq g(b)-g(a)+h(b)-h(a).
}
Upon taking the supremum with respect to all partitions $\tau$ of $[a,b]$, it follows that
\Eq{*}{
  V^{2\max(\Phi,\Psi)}_{[a,b]}f\leq g(b)-g(a)+h(b)-h(a)<\infty.
}
Hence $f$ is of bounded $2\max(\Phi,\Psi)$-total variation on $[a,b]$.
\end{proof}

The particular case $\Phi=\Psi$ of the above result yields the following statement.

\Cor{JDT1}{Let $\Phi\in\E(I)$. If $f:I\to\R$ is a delta-$\Phi$-monotone function, then the total $2\Phi$-variation of $f$ is finite on every compact subinterval of $I$.}

\Thm{JDT2}{Let $\Phi\in\E(I)$ let $f:I\to\R$ such that the total $2\Phi$-variation of $f$ on is finite on every compact subinterval of $I$. Then, for all $a\in I$, $f$ is a delta-$\Phi$-monotone function on $I\cap\,]a,\infty[$.}

\begin{proof}
Assume that the total $2\Phi$-variation of $f$ on every compact subinterval of $I$
is finite.

Let $a\in I$ be an arbitrarily fixed point and, for $x\in I$, $x>a$, define 
\Eq{*}{
  g(x):=\frac12\big(V^{2\Phi}_{[a,x]}f+f(x)\big) \qquad\mbox{and}\qquad
  h(x):=\frac12\big(V^{2\Phi}_{[a,x]}f-f(x)\big).
}
Then, we immediately have that $f=g-h$. 

Then, based on the \lem{JDT}, for $a<x<y$, we get
\Eq{*}{
  V^{2\Phi}_{[a,x]}f+f(x)-f(y)-2\Phi(y-x)
  &\leq V^{2\Phi}_{[a,x]}f+|f(x)-f(y)|-2\Phi(y-x)\\
  &\leq V^{2\Phi}_{[a,x]}f+V^{2\Phi}_{[x,y]}f\leq V^{2\Phi}_{[a,y]}f.
}
Rearranging this inequality, it follows that
\Eq{*}{
  g(x)\leq g(y)+\Phi(y-x),
}
which proves that $g$ is $\Phi$-monotone. Similarly, we can see that $h$ is 
also $\Phi$-monotone. This, together with the identity $f=g-h$ show that $f$ is
delta-$\Phi$-monotone function on $I\cap\,]a,\infty[$.
\end{proof}

\section{Individual error functions}

In this section we shall characterize the elements of the classes $\M(I)$ and $\H(I)$. For this purpose, given a function $f:I\to\R$, we define the following two extended real valued error functions on $[0,\ell(I)[\,$:
\Eq{*}{
\Phi_f^\sigma(u):=\sup_{x\in I\cap(I-u)}(f(x)-f(x+u))_+\qquad\mbox{and}\qquad
\Phi_f^\alpha(u):=\sup_{x\in I\cap(I-u)}|f(x)-f(x+u)|.
}
Here, the positive part of a real number $c$ is defined as $c_+:=\max(c,0)$.

\Thm{T}{Let $f:I\to\R$. Then we have the following two statements.
\begin{enumerate}[(i)]
 \item $f\in\M(I)$ if and only if $\Phi_f^\sigma$ is finite valued on $[0,\ell(I)[\,$. Additionally, if for some $\Phi\in\E(I)$, we have $f\in\M_\Phi(I)$, then $\Phi_f^\sigma\leq \Phi$.
 \item $f\in\H(I)$ if and only if $\Phi_f^\alpha$ is finite valued on $[0,\ell(I)[\,$. Additionally, if for some $\Phi\in\E(I)$, we have $f\in\H_\Phi(I)$, then $\Phi_f^\alpha\leq \Phi$.
\end{enumerate}
}

\begin{proof} Assume that $f\in\M(I)$. Then there exists an error function $\Phi\in\E(I)$ such that $f\in\M_\Phi(I)$. Therefore, for all $u\in[0,\ell(I)[\,$, 
\Eq{*}{
  f(x)-f(x+u)\leq \Phi(u) \qquad(x\in I\cap(I-u)).
}
Using that $\Phi$ is nonnegative, we get
\Eq{*}{
  (f(x)-f(x+u))_+=\max(f(x)-f(x+u),0)\leq \Phi(u) \qquad(x\in I\cap(I-u)).
}
Upon taking the supremum with respect to $x\in I\cap(I-u)$, we get that $\Phi_f^\sigma(u)\leq \Phi(u)$, which proves that $\Phi_f^\sigma$ has finite values.

Conversely, assume that $\Phi_f^\sigma$ has finite values. Then $\Phi_f^\sigma\in\E(I)$ and, for all $u\in[0,\ell(I)[\,$ and $x\in I\cap(I-u)$,
\Eq{*}{
  f(x)-f(x+u)\leq \Phi_f^\sigma(u).
}
This shows that $f$ is $\Phi_f^\sigma$-monotone, and hence, $f\in\M(I)$.

The proof about the second assertion is very similar and therefore, it is omitted.
\end{proof}

\Thm{PEF}{Let $f:I\to\R$. Then we have the following two statements. 
\begin{enumerate}[(i)]
 \item If $f\in\M(I)$, then $\Phi_{f}^\sigma$ is subadditive. Additionally, if for some $\Phi\in\E(I)$, we have $f\in\M_\Phi(I)$, then $\Phi_{f}^\sigma\leq\Phi^\sigma$.
 \item Provided that $I$ is unbounded, if $f\in\H(I)$, then $\Phi_{f}^\alpha$ is absolutely subadditive. Additionally, if for some $\Phi\in\E(I)$, $f\in\H_\Phi(I)$, then $\Phi_{f}^\alpha\leq\Phi^\alpha$.
\end{enumerate}
}

\begin{proof}
Let $f\in\M(I)$ and let $u,v\in\R_+$ such that $u+v<\ell(I)$. By the definition of $\Phi_f^\sigma$ and the subadditivity of the function $(\cdot)_+$, we have
\Eq{*}
{\Phi_f^\sigma(u+v)
&=\sup_{x\in I\cap(I-(u+v))}\big(f(x)-f(x+(u+v))\big)_+\\
&\leq\sup_{x\in I\cap(I-(u+v))}\big(f(x)-f(x+u)\big)_+
   +\big(f(x+u)-f(x+(u+v))\big)_+\\
&\leq \sup_{x\in I\cap(I-(u+v))}\big(f(x)-f(x+u)\big)_+
   +\sup_{x\in I\cap(I-(u+v))}\big(f(x+u)-f(x+(u+v))\big)_+\\
&\leq \sup_{x\in I\cap(I-u)}\big(f(x)-f(x+u)\big)_+
   +\sup_{x\in (I-u)\cap(I-(u+v))}\big(f(x+u)-f(x+(u+v))\big)_+\\
&=\Phi_f^\sigma(u)+\Phi_f^\sigma(v),
}
which establishes the subadditivity of $\Phi_f^\sigma$.

Since,$\Phi^\sigma$ is the largest subadditive function satisfying the inequality 
$\Phi^\sigma\leq\Phi$. Therefore, the inequality $\Phi_{f}^\sigma\leq\Phi^\sigma$ follows.

Suppose that the interval $I$ is unbounded from above. Let $f\in\H(I)$ and let $u,v\in\R$ such that $|u|,|v|,|u+v|<\ell(I)$. We may assume that $u+v$ is nonnegative (otherwise we replace $u$ by $(-u)$ and $v$ by $(-v)$). Therefore, by the unboundedness of $I$, we have that $I\subseteq I-(u+v)$ holds. Then at least one of the valuse $u$ and $v$ is nonnegative. By symmetry, we may also assume that $u\geq0$ and thus we have $I\subseteq I-u$. In this case, $x+u\in I$ for all $x\in I$. By the definition of $\Phi_f^\sigma$ and the subadditivity of the function $|\cdot|$, we have
\Eq{*}
{\Phi_f^\sigma(u+v)
&=\sup_{x\in I\cap(I-(u+v))}\big|f(x)-f(x+(u+v))\big|
 =\sup_{x\in I}\big|f(x)-f(x+(u+v))\big|\\
&\leq\sup_{x\in I}\big|f(x)-f(x+u)\big|
   +\big|f(x+u)-f(x+(u+v))\big|\\
&\leq \sup_{x\in I}\big|f(x)-f(x+u)\big|
   +\sup_{x\in I}\big|f(x+u)-f(x+(u+v))\big|\\
&\leq \sup_{x\in I\cap(I-u)}\big|f(x)-f(x+u)\big|
   +\sup_{x\in (I-u)\cap(I-(u+v))}\big|f(x+u)-f(x+(u+v))\big|\\
&=\Phi_f^\sigma(|u|)+\Phi_f^\sigma(|v|),
}
which establishes the absolute subadditivity of $\Phi_f^\alpha$ if $I$ is unbounded from above. The argument for the remaining case is similar, therefore it is left to the reader.

Since, $\Phi^\alpha$ is the largest absolutely subadditive function satisfying the inequality $\Phi^\alpha\leq\Phi$. Therefore, the inequality $\Phi_f^\alpha\leq\Phi^\alpha$ follows.
\end{proof}


\end{document}